\documentclass[letterpaper, 10 pt, conference]{ieeeconf}
\IEEEoverridecommandlockouts                              
\let\labelindent\relax

\usepackage{amsthm}
\usepackage[utf8]{inputenc}
\usepackage{amsmath}
\usepackage{amssymb}
\usepackage{xcolor}
\usepackage{enumitem}
\usepackage{comment}
\usepackage{graphicx} 
\usepackage{fixltx2e}
\usepackage{caption}
\usepackage{subcaption}
\usepackage{breqn}
\usepackage{stfloats}
\usepackage{hyperref}
\usepackage{xurl}

\newcommand{\B}[1]{\if#1\relax\bm{#1}\else\mathbf{#1}\fi} 
\newcommand{\R}{\mathbb{R}}
\newcommand{\I}{\mathbb{I}}

\newcommand{\sN}{\mathcal{N}}
\newcommand{\sKL}{\mathcal{KL}}
\newcommand{\sK}{\mathcal{K}}
\newcommand{\abs}[1]{\vert #1 \vert}
\newcommand{\norm}[1]{\Vert #1 \Vert}

\newcommand{\T}{^{\mathsf{T}}}
\newcommand{\maj}[1]{\lceil #1 \rceil}
\newcommand{\mmaj}[1]{\lceil #1 \rceil_{\text{M}}}

\newtheorem{definition}{Definition}

\newtheorem{lemma}{Lemma}

\newtheorem{remark}{Remark}
\newtheorem{proposition}{Proposition}

\title{On the design of multiplex control to reject disturbances in nonlinear network systems affected by heterogeneous delays}

\author{Shihao Xie$^{1}$ and Giovanni Russo$^{2}$
\thanks{$^{1}$Shihao Xie is with School of Electrical and Electronic Engineering, University College Dublin, Ireland
        {\tt\small shihao.xie1@ucdconnect.ie}}%
\thanks{$^{2}$Giovanni Russo is with the Department of Information and Electrical Engineering and Applied Mathematics, University of Salerno, Italy
        {\tt\small giovarusso@unisa.it}}%
}

\begin{document}
\maketitle
\thispagestyle{empty}
\pagestyle{empty}

\begin{abstract}
We consider the problem of designing control protocols for nonlinear network systems affected by heterogeneous, time-varying delays and disturbances. For these networks, the goal is to reject polynomial disturbances affecting the agents and to guarantee the fulfilment of some desired network behaviour. To satisfy these requirements, we propose an integral control design implemented via a multiplex architecture. We give sufficient conditions for the desired disturbance rejection and stability properties by leveraging tools from contraction theory. We illustrate the effectiveness of the results via a numerical example that involves the control of a multi-terminal high-voltage DC grid.
\end{abstract}

\section{Introduction}
Over the past decades, the size and complexity of network systems have considerably evolved thanks to the rapid development of computing and communication technologies. Much research efforts have been devoted to study collective behaviours such as consensus, synchronisation, formation control. A key challenge when designing the control protocols is to achieve desired behaviours despite imperfect communications, exogenous disturbances and delays.

In this context, we study the problem of designing distributed integral control protocols that guarantee the fulfilment of the desired network behaviour, while rejecting certain classes of disturbances. These requirements are captured via an Input-to-State Stability (ISS) property and we give sufficient conditions for this property based on non-Euclidean contraction theory.



\paragraph*{Related works} the design of integral control protocols for network systems that are able to reject constant disturbances has been investigated in e.g. \cite{silva2021string, knorn2014passivity}. In \cite{lombana2016multiplex}, a PI controller is delivered via multiplex architecture to achieve consensus. Recently, in \cite{xie2022design}, integral actions delivered by multiplex architecture with multiple layers are shown to be effective in rejecting higher order polynomial disturbances while guaranteeing a {\em scalability} property. The results in this paper are based on ideas from contraction theory \cite{lohmiller1998contraction}, particularly leveraging the use of non-Euclidean norms \cite{jafarpour2022robust, davydov2022non, russo2010global}.
We refer to \cite{aminzare2014contraction, di2016convergence, tsukamoto2021contraction} and references therein for details. In the context of delay-free networks, leveraging contraction theory, conditions for the synthesis of distributed controls using separable control contraction metrics are given in \cite{shiromoto2018distributed}; contracting recurrent network is introduced in \cite{revay2021recurrent} with guarantees of stability and robustness. For network systems affected by delays, \cite{wang2006contraction} shows the preservation of contraction for a time-delayed network using Euclidean contraction metric and, in \cite{xie2022design}, conditions are given for networks with homogeneous delays. 


\paragraph*{Statement of Contributions} we present a distributed multiplex integral control design for nonlinear network systems affected by both heterogeneous time-varying delays and disturbances (possibly with polynomial components). The goal of the control protocol is to guarantee, for the network: (i) rejection of polynomial disturbances; (ii) the fulfilment of some desired behaviours. These properties are rigorously formalised in Section \ref{sec:problem_set-up}. Specifically, our technical contributions can be summarised as follows: (i) we formalise the control problem as an Input-{Output} Stability problem and give sufficient conditions to assess this property. While the results of this paper leverage some of the tools from \cite{xie2022design}, here, differently from \cite{xie2022design}, we consider a weaker stability property for networks with heterogeneous delays. {In \cite{xie2022design}, although a stronger scalability property was considered, results were obtained by assuming that the network has homogeneous delays}; (ii) we show that our results can serve as design guidelines for the control protocol; (iii) the results are validated on the problem of designing a control protocol for multi-terminal high-voltage DC (MTDC) grid and we show how the conditions can be fulfilled by solving a convex optimisation problem. Simulations confirm the effectiveness of the results.

\section{Mathematical preliminaries}
Given a norm $\norm{\cdot}$, we denote by $\norm{A}$ its induced matrix norm with respect to a $m \times m$ real matrix $A$ and the corresponding matrix measure $\mu(A)=\lim_{h\rightarrow 0^+}\frac{\norm{I+hA}-1}{h}$. The symmetric part of a matrix $A$ is denoted as $[A]_s:=\frac{A+A\T}{2}$. For $\eta\in\R_{>0}^n$, we define a diagonal matrix by $[\eta]\in\R^{n\times n}$ with $[\eta]_{ii}=\eta_i$, $i\in \{1,\ldots,n\}$. The diagonally weighted $\ell_\infty$-norm of $x\in\R^{n}$ is defined as $\norm{x}_{\infty,[\eta]^{-1}}:=\max_i\{\abs{x_i}/\eta_i\}$ with the induced matrix norm $\norm{A}_{\infty,[\eta]^{-1}}:=\max_i\{\sum_j\frac{\eta_j}{\eta_i}\abs{A_{ij}}\}$ and matrix measure $\mu_{\infty,[\eta]^{-1}}(A):=\max_i\{ A_{ii}+\sum_{j\ne i}\frac{\eta_j}{\eta_i}\abs{A_{ij}}\}$. We denote by $\I_n$ the $n\times n$ identity matrix, by $0_{m\times n}$ the $m\times n$ zero matrix (if $m=n$ we simply write $0_n$) and by $\mathbf{1}_n\in\R^n$ the one vector. Let $f$ be a sufficiently smooth function, we denote by $f^{(n)}$ the $n$-th derivative of $f$. We recall that a continuous function $\alpha: [0,a)\rightarrow [0,\infty)$ is said to belong to class $\sK$ if it is strictly increasing and $\alpha(0)=0$. It is said to belong to class $\sK_\infty$ if $a=\infty$ and $\alpha(r)\rightarrow \infty$ as $r \rightarrow \infty$. A continuous function $\beta: [0,a)\times [0,\infty)\rightarrow [0,\infty)$ is said to belong to class $\sKL$ if, for each fixed $s$, the mapping $\beta(r,s)$ belongs to class $\sK$ with respect to $r$ and, for each fixed $r$, the mapping $\beta(r,s)$ is decreasing with respect to $s$ and $\beta(r,s)\rightarrow 0$ as $s\rightarrow \infty$.


The following results, originally introduced in \cite{5717887}, can be found in its current form in \cite{FB-CTDS}.
\begin{lemma}\label{lem:matrix norm}
Given $r$ positive integers $n_1,\ldots, n_r$ such that $n_1+\cdots+n_r=n$. Consider the vector $x:=[x_1\T,\ldots,x_r\T]\T\in \R^n$, $x_i \in \mathbb{R}^{n_i}$. We let the composite norm $\norm{x}_{\text{cmpst}} := \norm{\left[ \norm{x_1}_{1},\ldots,\norm{x_r}_{r}\right]}_{\text{agg}}$, with $\norm{\cdot}_{i}$($\norm{\cdot}_{\text{agg}}$) being local(aggregating) norms on $\R^{n_i}$($\R^{n}$), and induced matrix norm $\norm{\cdot}_i$($\norm{\cdot}_{\text{agg}}$) and matrix measure $\mu_i(\cdot)$($\mu_{\text{agg}}(\cdot)$). Finally, given an $n\times n$ block matrix $A$ with $A_{ij}\in \R^{n_i\times n_j}$, define: 
\begin{enumerate}
    \item[(i)] the aggregate majorant $\maj{A}$ with $(\maj{A})_{ij}=\norm{A_{ij}}_{ij}$;
    \item[(ii)] the aggregate Metzler majorant $\mmaj{A}$ with
    \begin{equation*}
  (\mmaj{A})_{ij}=
    \begin{cases}
      \norm{A_{ij}}_{ij}, \text{if } j\ne i,\\
      \mu_i(A_{ii}), \text{if } j=i,
    \end{cases}       
\end{equation*}
\end{enumerate}
where $\norm{A_{ij}}_{{ij}}:=\sup_{\norm{x_j}_{j}=1}\norm{A_{ij}x_j}_{i}$. \\
\noindent Then, 
\begin{itemize}
\item[(i)] the composite norm $\norm{\cdot}_{\text{cmpst}}$ is well-defined, i.e. satisfying the norm properties;
\item[(ii)] If the aggregating norm $\norm{\cdot}_{\text{agg}}$ is monotonic, then: 
\begin{enumerate}
    \item $\mu_{\text{cmpst}}(A) \le \mu_{\text{agg}}(\mmaj{A})$;
    \item $\norm{A}_{\text{cmpst}} \le \norm{\maj{A}}_{\text{agg}}$.
\end{enumerate}
\end{itemize}

\end{lemma}
If the norms $\norm{\cdot}_i$, $\norm{\cdot}_{\text{agg}}$ in Lemma \ref{lem:matrix norm} are $p$-norms (with the same $p$) then $\norm{\cdot}_{\text{cmpst}}$ is again a $p$-norm defined on a larger space. The next lemma follows from \cite[Theorem $2.4$]{wen2008generalized}.
\begin{lemma}\label{lem: halanay}
Let $u:[-\tau_0,+\infty)\rightarrow\R_{\ge 0}$ , $\tau_0<+\infty$ and assume that
$$
D^+u(t) \le au(t)+b \sup_{t-\tau(t) \le s \le t}u(s)+c, \ \ t\ge t_0 
$$
with: (i) $\tau(t)$ being bounded and non-negative, i.e. $0{\le}\tau(t)\le\tau_0$, $\forall t$; (ii) $u(t)=\abs{\varphi(t)}$, $\forall t\in [t_0-\tau_0,t_0]$ where $\varphi(t)$ is bounded {and continuous} in $[t_0-\tau_0,t_0]$; (iii) $a < 0$, $b \ge 0$ and $c \ge 0$. Assume that there exists some $\sigma>0$ such that $a+b \le -\sigma <0, \forall t\ge t_0$. Then:
$$
u(t) \le \sup_{t_0-\tau_0 \le s \le t_0}u(s)e^{-\hat \lambda t}+\frac{c}{\sigma}
$$
where $\hat \lambda:=\inf_{t\ge t_0}\{\lambda|\lambda(t) + a+be^{\lambda(t)\tau(t)}=0\}$ is positive.
\end{lemma}

\section{System description and problem formulation}\label{sec:problem_set-up}
\begin{figure}
\centering
\includegraphics[width=0.95\columnwidth]{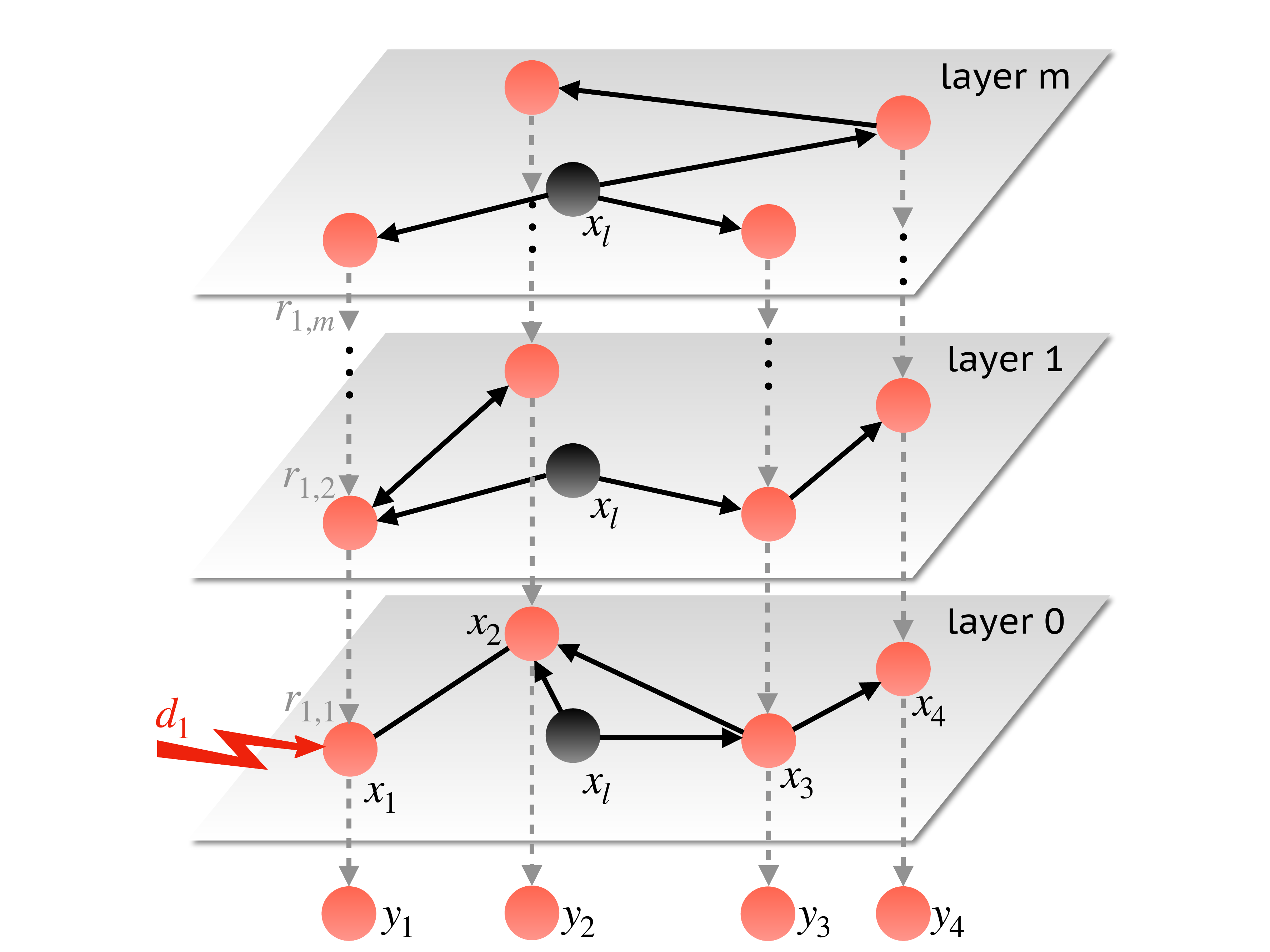}
\caption{An example of the multiplex architecture (layer $0$ to layer $m$). The agents are represented by orange nodes and the leader is represented by the black node. {Layers are not required to have the same topology. Each layer outputs a signal to the layer immediately below (if any), i.e. layer $k$ outputs $[r_{1,k}\T,\ldots,r_{N,k}\T]\T$ to layer $k-1$.} Only one disturbance, e.g. $d_1$ on agent $1$, is shown.}
\label{fig: structure}
\end{figure}

Consider a network system comprised of $N$ agents with the dynamics of the $i$-th agent described by
\begin{align}\label{equ: dynamics}
\begin{split}
    &\dot{x}_i(t)=f_i(x_i(t),t)+u_i(t)+d_i(t) \\
    &{y_i(t)=g_i(x_i(t))}
\end{split}
\end{align}
$t\ge t_0$, $i\in\{1,\ldots,N\}$. In the above expression, $x_i(t)\in \R^n$ denotes the agent state, $u_i(t)\in \R^n$ denotes the control input, $f_i: \R^n \times \R \rightarrow \R^n$ is a smooth function and {$g_i: \R^n\rightarrow \R^r$ is the output function, which we assume to be Lipschitz}. The term $d_i(t)\in \R^n$ models the external disturbance of the form: 
\begin{align}\label{equ: disturbance}
    d_i(t)=w_i(t) + \bar{d}_i(t) := w_i(t)+\sum_{k=0}^{m-1} \bar 
    d_{i,k}\cdot t^k
\end{align}
where $\bar{d}_i(t)$ represents the polynomial component of the disturbance of the order $m-1$ ($m\in \mathbb{Z}_{> 0}$) with $\bar d_{i,k}$ being constant vectors and $w_i(t)$ is a piecewise continuous signal capturing \emph{residual} terms in the disturbance that are not modelled with the polynomial. In order to reject polynomials up to order $m-1$, we design the control protocol $u_i(t)$ following \cite{xie2022design} with integral actions delivered by $m$ multiplex layers as illustrated in Figure \ref{fig: structure}. 
\begin{align}\label{equ: control}
\begin{split}
    u_i(t)&=h_{i,0}(x_i,\{x_{j}\}_{j\in \sN_i},x_l,t)+h_{i,0}^{(\tau)}(x_i,\{x_{j}\}_{j\in \sN_i},x_l,t)\\
    &+ r_{i,1}(t)\\
    \dot{r}_{i,1}(t)&=h_{i,1}(x_i,\{x_j\}_{j\in \sN_i},x_l,t)+h_{i,1}^{(\tau)}(x_i,\{x_j\}_{j\in \sN_i},x_l,t)\\
    &+ r_{i,2}(t)\\
    \vdots\\
    \dot{r}_{i,m}(t)&=h_{i,m}(x_i,\{x_j\}_{j\in \sN_i},x_l,t)+h_{i,m}^{(\tau)}(x_i,\{x_j\}_{j\in \sN_i},x_l,t)
\end{split}
\end{align}
where $\{x_j\}_{j\in\sN_i}$ denotes the state of the neighbours of agent $i$, $x_l(t)$ is the {exhogenous reference signal} from e.g., a leader. The functions  $h_{i,k}, h_{i,k}^{(\tau)}, \forall k$ are smooth coupling functions on $k$-th multiplex layer modelling delay-free and delayed communications from neighbours of agent $i$ and the leader. For delayed communications, we consider time-varying delays $\tau_{ij}(t), \tau_{il}(t)\in (0,\tau_{\max}]$ when information is transmitted to agent $i$ from agent $j$ and from the leader, respectively. Note that in general,  $\tau_{ij}(t)\ne \tau_{ji}(t)$. In what follows, time dependence inside these coupling functions are omitted for notational convenience. 
\begin{remark}\label{rmk: delay}
Protocol in \eqref{equ: control} {can be used for both leaderless and leader-follower networks.} It arises in a wide range of applications. For example, the classic diffusive-type protocol can be written as in \eqref{equ: control} with $u_i(t)=h_{i,0}^{(\tau)}(x_i,\{x_{j}\}_{j\in \sN_i},x_l,t)=\sum_{j\in \sN_i}a_{ij}(x_j(t-\tau_{ij})-x_i(t-\tau_{ij}))$, i.e. when all the communications are delayed and no integral actions are applied {in a leaderless network.}
\end{remark}
Let $x(t)=[x_1\T(t),\ldots,x_N\T(t)]\T$, $u(t)=[u_1\T(t),\ldots,u_N\T(t)]\T$, {$y(t)=[y_1\T(t),\ldots,y_N\T(t)]\T$} and $d(t)=w(t)+\bar{d}(t)$ where $w(t)=[w_1\T(t),\ldots,w_N\T(t)]\T$, $\bar{d}(t)=[\bar{d}_1\T(t),\ldots,\bar{d}_N\T(t)]\T$, the interconnected system \eqref{equ: dynamics} can then be written in a compact form as
\begin{align}\label{equ: system_dynamics}
\begin{split}
    &\dot x(t)=f(x(t),t)+u(t)+d(t)\\
    & {y(t)=g(x(t))} 
\end{split}
\end{align}
$t\ge t_0$, where $f(x(t),t)=[f_1\T(x_1(t),t), \ldots, f_N\T(x_N(t),t)]\T$ and {$g(x(t))=[g_1\T(x_1(t)), \ldots, g_N\T(x_N(t))]\T$}. We also let $r_k(t)=[r_{1,k}\T(t),\ldots,r_{N,k}\T(t)]\T$, $ k\in\{1,\ldots,m\}$. Next, we state the control goal in terms of the desired solution of the network when there are no disturbances and delays. {The desired solution 
is the solution of the unperturbed system satisfying: (i) $\dot{x}_i^{\ast}=f_i(x_i^\ast,t)$, $\forall i$, {$\forall t\ge t_0$} and (ii) $r_{i,k}^\ast(t)=0$, $\forall i, \forall k$. The desired output is then $y^\ast(t):=[y_1^{*\mathsf{T}}(t), \dots,y_N^{\ast\mathsf{T}}(t)]\T$, with ${y}^\ast_i(t)=g_i(x_i^\ast(t)), \forall i$.}
In what follows, we simply term $x^\ast(t):=[x_1^{\ast\mathsf{T}}(t),\ldots,x_N^{\ast\mathsf{T}}(t)]\T$ as desired solution. We are ready to give the following:
\begin{definition}\label{def: L_inf}
The closed-loop system \eqref{equ: system_dynamics} affected by disturbance $d(t)=w(t) + \bar{d}(t)$ is Input-{Output} Stable with respect to $w(t)$ if there exists class $\sKL$ functions $\alpha(\cdot,\cdot)$, $\beta(\cdot,\cdot)$, a class $\sK$ function $\gamma(\cdot)$, such that 
    \begin{align*}
        \begin{split}
            &{\norm{y(t)-y^*(t)}} \le\alpha\left(\sup_{t_0-\tau_{\max}\le s \le t_0}\norm{x(s)-x^*(s)},t-t_0\right)\\
            &+\beta\left(\sup_{t_0-\tau_{\max}\le s \le t_0}\sum_{k=1}^m\norm{r_{k}(s)+\bar{d}^{(k-1)}(s)},t-t_0\right)\\
            &+\gamma\left(\sup_t\norm{w(t)}\right)
        \end{split}
    \end{align*}
\end{definition}
\noindent holds $\forall t \ge t_0$, $\forall x(s)$, $\forall r_k(s), k\in\{1,\ldots,m\}$, where $x(s)=\varphi(s), r_{k}(s)=\phi_{k}(s)$ with $\varphi(s), \phi_{k}(s)$ being continuous and bounded functions in $[t_0-\tau_{\max},t_0]$.

\section{Technical results}\label{sec: technical}
In this section, we give a sufficient condition guaranteeing that the closed-loop system (\ref{equ: system_dynamics}) affected by disturbances of the form (\ref{equ: disturbance}) is Input-Output Stable. The results are stated in terms of:
\begin{enumerate}
    \item a composite norm (\cite[Section 2.4.4]{FB-CTDS}) $\norm{x}_{\text{cmpst}}=\norm{\left[ \norm{x_1}_{p},\ldots,\norm{x_N}_{p}\right]}_{\infty,[\eta]^{-1}}$ where $[\eta]^{-1} \in \R^{N\times N}$ is a diagonal weighting matrix with $\eta\in\R^N_{>0}$;
    \item a block diagonal coordinate transformation matrix $T=\text{diag}\{T_1,\ldots,T_N\}$. 
\end{enumerate}
For the statement of our result, it is also useful to {\em relabel} the delays affecting the network. Specifically, we define $\tau_k(t), k = 1, \ldots, q$ with $q \le N^2$, as an element of the set 
$$\{\tau_{ij}(t): i, j = 1, \ldots, N, i \ne j\}\bigcup\{\tau_{il}(t):i = 1, \ldots, N\}$$
\begin{proposition}\label{prop: stable}
Consider the closed-loop system (\ref{equ: dynamics}) affected by external disturbance (\ref{equ: disturbance}). Assume that, $\forall t\ge t_0$ and for some matrices $T_i, T_j, i,j =1,\ldots, N$, the following conditions are satisfied for some $0<\underline{\sigma}< \bar{\sigma}<\infty$: 
\begin{itemize}
    \item[C1] $h_{i,k}(x_i^*,\{x_{j}^\ast\}_{j\in \sN_i},x_l,t)=h_{i,k}^{(\tau)}(x_i^*,\{x_{j}^\ast\}_{j\in \sN_i},x_l,t)=0$, $\forall i, k$;
    \item[C2] $\mu_p(T_i\tilde A_{ii}(t)T_i^{-1})+\sum_{j}\frac{\eta_j}{\eta_i}\norm{T_i\tilde A_{ij}(t)T_j^{-1}}_p\le -\bar{\sigma}$, $\forall i$, {$\forall x$, $\forall x_l$};
    \item[C3] $\sum_{k=1}^q\sum_{j}\frac{\eta_j}{\eta_i}\norm{T_i(\tilde{B}_{k}(t))_{ij}T_j^{-1}}_p\le \underline{\sigma}$, $\forall i$, {$\forall x$, $\forall x_l$}.
\end{itemize}
In the above expression, 
\begin{align}\label{equ: matrices}
\begin{split}
    \tilde{A}_{ii}(t) & =\left[\begin{matrix} \frac{\partial f_i}{\partial x_i}+\frac{\partial 				h_{i,0}}{\partial x_i} & \I_n & 0_n & \cdots & 0_n\\
     			\frac{\partial h_{i,1}}{\partial x_i} & 0_n & \I_n & \cdots & 0_n\\ 						\vdots & \vdots & \vdots & \ddots & \vdots\\
     			\frac{\partial h_{i,m-1}}{\partial x_i} & 0_n & 0_n & \cdots & \I_n\\
      			\frac{\partial h_{i,m}}{\partial x_i} & 0_n & 0_n & \cdots & 0_n 	\end{matrix}\right]\\
      			\tilde{A}_{ij}(t) & =\left[\begin{matrix} \frac{\partial h_{i,0}}{\partial x_j} & 0_n & \cdots & 0_n\\
      			\vdots & \vdots & \ddots & \vdots\\
      			\frac{\partial h_{i,m}}{\partial x_j} & 0_n & \cdots & 0_n 	\end{matrix}\right]  	\\
(\tilde{B}_k)_{ij}(t)& =\left[\begin{matrix} \frac{\partial h_{i,0}^{(\tau)}}{\partial x_j} & 		0_n & \cdots & 0_n\\
      			\vdots & \vdots & \ddots & \vdots\\
      			\frac{\partial h_{i,m}^{(\tau)}}{\partial x_j} & 0_n & \cdots & 0_n  \end{matrix}\right]
\end{split}
\end{align}

Then, the system is Input-Output Stable. In particular,
\begin{align*}
\begin{split}
            &{\norm{y(t)-y^*(t)}_{\text{cmpst}}} \le \norm{T^{-1}}_{\text{cmpst}}\norm{T}_{\text{cmpst}}\cdot\\
            &\cdot e^{-\lambda (t-t_0)}\biggl(\sup_{t_0-\tau_{\max}\le s \le t_0}\norm{x(s)-x^*(s)}_{\text{cmpst}}\\
            &+\sup_{t_0-\tau_{\max}\le s \le t_0}\sum_{k=1}^m\Vert \sum_{b=0}^{m-k} \frac{(m-1-b)!}{(m-k-b)!}\cdot \bar d_{m-1-b}\cdot s^{m-k-b} \\
            &+ r_{k}(s)  \Vert_{\text{cmpst}}\biggr)+\frac{\norm{T^{-1}}_{\text{cmpst}}\norm{T}_{\text{cmpst}}}{\bar{\sigma}-\underline{\sigma}}\sup_t\norm{w(t)}_{\text{cmpst}}
\end{split}
\end{align*}
where {$\lambda>0$} {is the solution of the following equation:}
\begin{align}\label{equ: delay}
    \lambda-\bar{\sigma}+\underline{\sigma}e^{\lambda\tau_{\max}}=0
\end{align}
\end{proposition}
\begin{remark}\label{rmk: conditions}
Condition $C1$ implies that $u_i(t)=0$ at the desired solution which guarantees $x^\ast(t)$ is a solution of the unperturbed dynamics. $C2$ and $C3$ give conditions on the Jacobian of delay-free and delayed part of agent dynamics, respectively. These conditions are concerned with the couplings between agents and depend on the number of neighbours of the agents and on the coupling strength.
\end{remark}
\begin{remark}\label{rmk: desired_solution}
If $C2$ and $C3$ are satisfied, then the network is connective stable \cite[Chapter 2.1]{siljak2011decentralized}. In the special case where the agents are contracting and there are no delays, the conditions yield a known property of contracting systems: convergence of the network is achieved even without coupling (however, in this case, there is no guarantee that the solutions towards which the agents converge is the desired solution). 
\end{remark}
\begin{remark}
The diagonal weighting matrix $[\eta]^{-1}$ in $\norm{x}_{\text{cmpst}}$ can be carefully designed to achieve sharper bounds for the induced matrix norm, according to \cite{FB-CTDS}. Hence, we could find such matrix $[\eta]^{-1}$ to minimise the upper bound in $C2$ or $C3$. On the other hand, the selection of the coordinate transformation matrix $T$ is also of great significance as in many cases it turns out to be difficult to design controllers fulfilling conditions proposed for the original system. The conditions in Proposition \ref{prop: stable} not only provides guideline for the design of the controller, but also for the computation of matrix $T$, see also \cite{xie2022design}.
\end{remark}

\begin{remark}
Following Proposition \ref{prop: stable}, the closed-loop network has a convergence rate $\lambda$ depending on $\tau_{\max}$. As is highlighted in \eqref{equ: delay}, larger $\tau_{\max}$ yields lower $\lambda$ which means slower convergence towards the system desired solution. Hence \eqref{equ: delay} gives an implicit condition on communication delays for networks with desired convergence rates.
\end{remark}
{The proof is inspired by the derivations of \cite[Proposition $1$]{xie2022design}. Hence, some technical steps are omitted here for brevity.}

\begin{proof} Consider the augmented state $z_i(t)=[x_i\T(t), \zeta_{i,1}\T(t),\ldots,\zeta_{i,m}\T(t)]\T$ where
\begin{align*}
    \zeta_{i,k}(t)&=r_{i,k}(t) + \sum_{b=0}^{m-k} \frac{(m-1-b)!}{(m-k-b)!}\cdot \bar d_{i,m-1-b}\cdot t^{m-k-b}
\end{align*}
$k\in\{1,\cdots,m\}$. The desired solution is $z_i^\ast(t)=[x_i^{\ast \mathsf{T}}(t), 0_{1\times n},\ldots,0_{1\times n}]\T$. Hence the error dynamics $e_i(t):=z_i(t)-z^\ast_i(t)$ follows 
\begin{align*}
\begin{split}
    \dot{e}_i(t)&=\tilde{f}_i(z_i,t)-\tilde{f}_i(z_i^\ast,t)
    +\tilde h_i(z,t)+\tilde h_i^{(\tau)}(z,t)+ \tilde w_i(t)
\end{split}
\end{align*}
where $\tilde{f}_i(z_i,t)=[f_i\T(x_i,t), 0_{1\times n}, \ldots, 0_{1\times n}]\T$ and
\begin{align*}
    \begin{split}
        \tilde h_i(z,t) = \left[\begin{array}{*{20}c}
h_{i,0}(x_i,\{x_{j}\}_{j\in \sN_i},x_l,t) + \zeta_{i,1}(t) \\
h_{i,1}(x_i,\{x_{j}\}_{j\in \sN_i},x_l,t) + \zeta_{i,2}(t)\\
\vdots \\
h_{i,m-1}(x_i,\{x_{j}\}_{j\in \sN_i},x_l,t) + \zeta_{i,m}(t)\\
h_{i,m}(x_i,\{x_{j}\}_{j\in \sN_i},x_l,t) \end{array}\right]
         \end{split}
\end{align*}  
and
\begin{align*}
    \begin{split}
        \tilde h_i^{(\tau)}(z,t) = \left[\begin{array}{*{20}c}
h_{i,0}^{(\tau)}(x_i,\{x_{j}\}_{j\in \sN_i},x_l,t)\\
h_{i,1}^{(\tau)}(x_i,\{x_{j}\}_{j\in \sN_i},x_l,t)\\
\vdots \\
h_{i,m}^{(\tau)}(x_i,\{x_{j}\}_{j\in \sN_i},x_l,t)\end{array}\right]
    \end{split}
\end{align*}
and $\tilde w_i(t)=[w_i\T(t),0_{1\times n},\ldots,0_{1\times n}]\T$. 
Then following \cite{desoer1972measure}, let $\eta_i(\rho)=\rho z_i+(1-\rho)z_i^\ast$, $\eta(\rho)=[\eta_1\T(\rho),\ldots,\eta_N\T(\rho)]\T$ and we can rewrite the error dynamics in a compact form as 
\begin{align}\label{equ: scalable_error}
    \dot{e}(t)=A(t)e(t)+\sum_{k=1}^{q}B_k(t)e(t-\tau_k(t))+\tilde w(t)
\end{align}
where $\tilde w=[\tilde w_1\T(t),\ldots,\tilde w_N\T(t)]\T$. The Jacobian matrix $A(t)$ has entries $A_{ij}(t)=\int_0^1 \tilde A_{ij}(t) d\rho$ where $\tilde A_{ij}(t)= J_{\tilde{f}_i}(\eta_j(\rho),t)+J_{\tilde{h}_i}(\eta_j(\rho),t)$ and $B_k(t)$ with entries $(B_k)_{ij}(t)=\int_0^1 (\tilde B_k)_{ij}(t)d\rho$ where $(\tilde B_k)_{ij}(t)=J_{\tilde{h}_i^{(\tau)}}(\eta_j(\rho),t)$. Then let $\tilde{e}(t)=Te(t)$ where $T=\text{diag}\{T_1,\ldots,T_N\}$,
we have
\begin{align*}
    \dot{\tilde e}(t)=TA(t)T^{-1}\tilde e(t)+\sum_{k=1}^{q}TB_k(t)T^{-1}\tilde e(t-\tau_k(t))+T\tilde w(t)
\end{align*}
Taking the Dini derivative of $\norm{\tilde e(t)}_{\text{cmpst}}$, we obtain
\begin{align*}
    \begin{split}
        &D^+\norm{\tilde e(t)}_{\text{cmpst}}=\limsup_{h\rightarrow 0^+} \frac{1}{h}\left(\norm{\tilde e(t+h)}_{\text{cmpst}}-\norm{\tilde e(t)}_{\text{cmpst}}\right)\\
        &\le \mu_{\text{cmpst}}(TA(t)T^{-1} )\norm{\tilde e(t)}_{\text{cmpst}} +\sup_t\norm{T\tilde w(t)}_{\text{cmpst}} \\
        &\phantom{=}+\sum_{k=1}^q\norm{TB_k(t)T^{-1} }_{\text{cmpst}} \sup_{t-\tau_{\max}\le s \le t}\norm{\tilde e(s)}_{\text{cmpst}}
    \end{split}
\end{align*}
Condition $C2$ and $C3$ implies that $\mu_{\text{cmpst}}(TA(t)T^{-1})\le -\bar{\sigma}$ and $\sum_{k=1}^q\norm{TB_k(t)T^{-1} }_{\text{cmpst}}\le \underline{\sigma}$. {Then, following similar steps to the ones of the proof of} Proposition $1$ in \cite{xie2022design} and {leveraging the Lipschitz assumption of the output function},
we get the upper bound of the output deviation of \eqref{equ: dynamics} from its desired output:
\begin{align*}
        \begin{split}
            &{\norm{y(t)-y^*(t)}_{\text{cmpst}}} \le \\
            &\norm{T^{-1}}_{\text{cmpst}}\norm{T}_{\text{cmpst}}e^{-\lambda (t-t_0)}\biggl(\sup_{t_0-\tau_{\max}\le s \le t_0}\norm{x(s)-x^*(s)}_{\text{cmpst}}\\
            &+\sup_{t_0-\tau_{\max}\le s \le t_0}\sum_{k=1}^m\Vert \sum_{b=0}^{m-k} \frac{(m-1-b)!}{(m-k-b)!}\cdot \bar d_{m-1-b}\cdot s^{m-k-b} \\
            &+ r_{k}(s)  \Vert_{\text{cmpst}}\biggr)+\frac{\norm{T^{-1}}_{\text{cmpst}}\norm{T}_{\text{cmpst}}}{\bar{\sigma}-\underline{\sigma}}\sup_t\norm{w(t)}_{\text{cmpst}}
        \end{split}
\end{align*}
\end{proof} 

\section{Application Example}\label{sec: application}
\begin{figure}
\centering
\includegraphics[width=0.95\columnwidth]{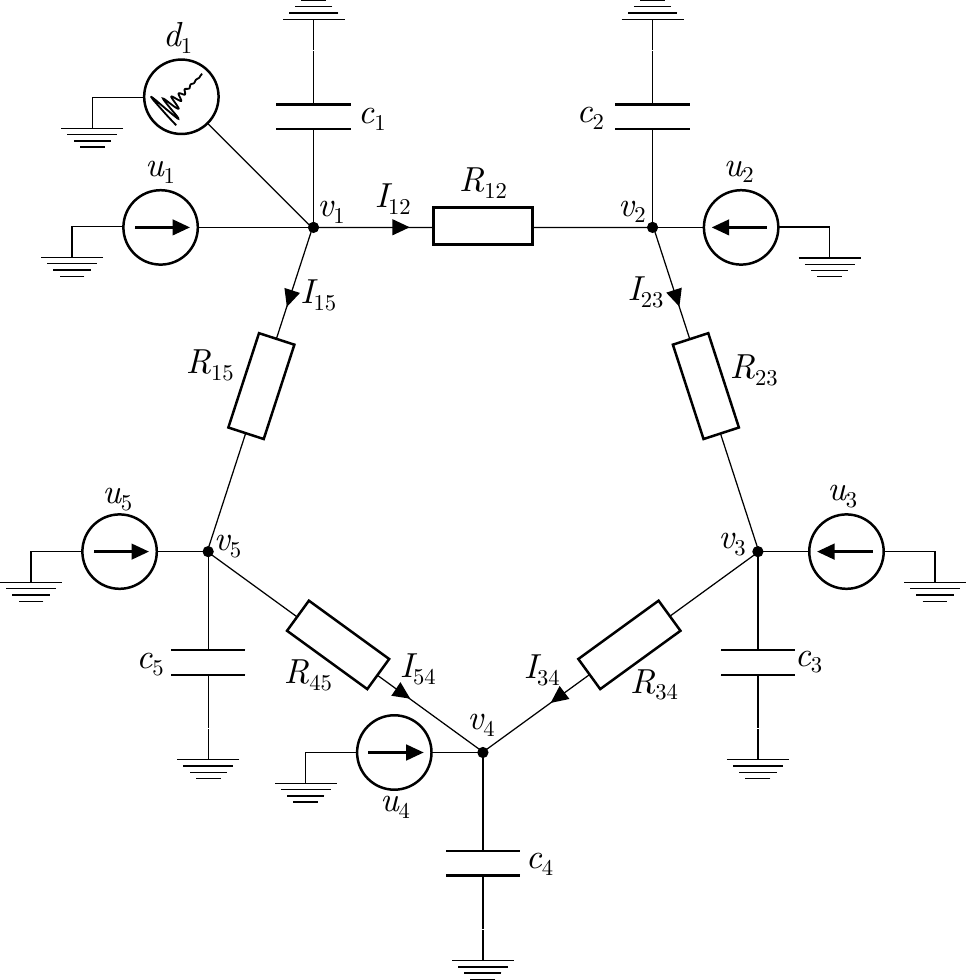}
\caption{Illustration of the MTDC composed of $5$ terminals. Terminal $1$ is affected by a disturbance $d_1$.}
\label{fig: topology}
\end{figure}

We consider an MTDC grid model from \cite{axelson2022online}. In our example, the grid has {$30$ terminals arranged in a ring (an example illustration of a ring consisting of $5$ terminals is given in Figure \ref{fig: topology}).} The dynamics of the terminals are described by
\begin{align}\label{equ: MTDC}
\begin{split}
    c_i\dot{v}_i(t)&=-\sum_{j\in\sN_i}I_{ij}(t)
    +u_i(t)+d_i(t), \ i\in\{1,\ldots,5\}\\
    y_i(t) &= v_i(t)
    \end{split}
\end{align}
where the current $I_{ij}(t)=\frac{1}{R_{ij}}(v_i(t)-v_j(t))$ according to Ohm's law and note that $v_j(t)$ is not delayed in this equation as communication is not required. In the above expression, $v_i$ denotes the voltage deviation of terminal $i$ from the nominal voltage $v^{\text{nom}}$ which is assumed to be identical for all terminals, $c_i$ is the capacitance of terminal $i$, $R_{ij}$ is the line resistance between terminal $i$ and terminal $j$ satisfying $R_{ij}=R_{ji}$, $u_i$ is the injected control current, $d_i$ is the disturbance due to e.g. load changes at the terminal and {$y_i$ is the output}. The design of control protocols to reject constant disturbances for \eqref{equ: MTDC} was considered in \cite{axelson2022graph} but without delays and in \cite{andreasson2014distributed} with constant delays. {In this example, the capacitance $c_i=1$mF and the resistance $R_{ij}=20\Omega$.} We now leverage Proposition \ref{prop: stable} to consider the case where delays are heterogeneous and the disturbances are polynomial. 
\subsection{Controller design}
In this example, we consider a first order disturbance, which can model e.g. the rapid increase of the current in the terminal caused by fault \cite{sneath2014fault}. To reject such disturbance before being diagnosed, we design the controller with $2$ multiplex layers following:
\begin{align}\label{equ: mtdc_control}
    \begin{split}
        u_i(t)&= -k_{0} v_i(t) -\sum_{j\in\sN_i}k_0^{(\tau)}(v_{i}(t-\tau_{ij}(t))-v_{j}(t-\tau_{ij}(t)))\\&+r_{i,1}(t)\\
        \dot{r}_{i,1}(t)&=-k_{1} v_i(t)  -\sum_{j\in\sN_i}k_1^{(\tau)}(v_{i}(t-\tau_{ij}(t))-v_{j}(t-\tau_{ij}(t)))\\&+ r_{i,2}(t)\\
        \dot{r}_{i,2}(t)&=-k_{2} v_i(t)-\sum_{j\in\sN_i}k_2^{(\tau)}(v_{i}(t-\tau_{ij}(t))-v_{j}(t-\tau_{ij}(t)))
    \end{split}
\end{align}
In the above expression, the delay occurs because the voltage information needs to be transmitted from terminal $j$ to terminal $i$ via communication (see also \cite{andreasson2014distributed}). For our design, we consider the composite norm $\norm{x}_{\text{cmpst}}=\norm{\left[ \norm{x_1}_{2},\ldots,\norm{x_N}_{2}\right]}_{\infty,{[\eta]^{-1}}}, {\eta\in\R_{>0}^N}$ and the coordinate transformation matrix $T$ with identical diagonal blocks
\begin{align*}
    T_i={\tilde{T}=}\left[\begin{matrix} 1 & \alpha & 0\\
                                        0 & 1 & \beta\\
                                        0 & 0 & 1\end{matrix}\right], \ \alpha, \beta\in\R, \ \forall i
\end{align*}
The desired solution of \eqref{equ: MTDC} is {$v_i^\ast(t)=0$, $\forall i$}. Hence, $C1$ is guaranteed by design of control \eqref{equ: mtdc_control}. We recast the fulfilment of $C2$ and $C3$ as an optimisation problem following
{
\begin{subequations}\label{equ: optimisation}
    \begin{align}
        &\min_{\xi} \ \mathcal{J}  \\
		&s.t. \quad k_{0}\ge0, k_{1}\ge 0, k_{2}\ge 0, k_{0}^{(\tau)}\ge 0, k_{1}^{(\tau)}\ge 0, k_{2}^{(\tau)}\ge 0 \\
		&\phantom{s.t. \quad} k_0+k_0^{(\tau)}>0, k_1+k_1^{(\tau)}>0, k_2+k_2^{(\tau)}>0 \\
		&\phantom{s.t. \quad} \bar{\sigma}>0, \underline{\sigma}\ge 0, \bar{\sigma}-\underline{\sigma}>0 \\
		&\phantom{s.t. \quad} \mu_2(\tilde T\tilde A_{ii}\tilde T^{-1})+\sum_{j}\frac{\eta_j}{\eta_i}\norm{\tilde T\tilde A_{ij}\tilde T^{-1}}_2\le -\bar{\sigma},\ \ \forall i,j\label{measure}\\
		&\phantom{s.t. \quad}  \sum_{k=1}^{q}\sum_{j}\frac{\eta_j}{\eta_i}\norm{\tilde T(\tilde{B}_{k})_{ij}\tilde T^{-1}}_2\le \underline{\sigma}, \ \ \forall i,j\label{norm}
    \end{align}
\end{subequations}
where $\xi:=[k_0, k_1, k_2, k_{0}^{(\tau)}, k_{1}^{(\tau)}, k_{2}^{(\tau)}, \bar{\sigma}, \underline{\sigma}]$ are the decision variables and $\mathcal{J}=-k_{0}^{(\tau)}-k_{1}^{(\tau)}-k_{2}^{(\tau)}$ is the cost function chosen, in this case, to maximise the coupling between terminals. In the above expression, the constraints on matrix measure and matrix norm ($\tilde{A}_{ii}, \tilde{A}_{ij}, (\tilde{B}_k)_{ij}$ are computed from \eqref{equ: matrices} and are constant matrices) can be recast as LMIs, for fixed $\alpha$, $\beta$ and $\eta$.
Indeed, \eqref{measure} is satisfied if 
\begin{align}\label{measure1}
\begin{split}
   & \mu_2(\tilde T\tilde A_{ii}\tilde T^{-1})\le -b_1, \ \ \forall i\\
   & \sum_{j}\frac{\eta_j}{\eta_i}\norm{\tilde T\tilde A_{ij}\tilde T^{-1}}_2\le b_2, \ \ \forall i,j\\
   & b_1> 0, b_2 \ge 0, -b_1+b_2\le -\bar{\sigma}
\end{split}
\end{align}
where $b_1, b_2$ are auxiliary variables. Also, note that due to the ring topology, for terminal $i$ there are two neighbours. It is then useful to recast $\sum_{j}\frac{\eta_j}{\eta_i}\norm{\tilde T\tilde A_{ij}\tilde T^{-1}}_2\le b_2$ as
\begin{align}\label{norm1}
\begin{split}
    & \frac{\eta_{j}}{\eta_i}\norm{\tilde T\tilde A_{i{j}}\tilde T^{-1}}_2  \le b_3, \ \ \forall i,j \\
    & b_3\ge 0, 2b_3 \le b_2
\end{split}
\end{align}
with auxiliary variable $b_3$. Analogously, we recast \eqref{norm} as
\begin{align}\label{norm2}
\begin{split}
  &  \norm{\tilde T(\tilde B_k)_{ii}\tilde T^{-1}}_2\le b_4, \ \ \forall i \\
   & \frac{\eta_{j}}{\eta_i}\norm{\tilde T(\tilde B_k)_{i{j}}\tilde T^{-1}}_2\le b_5, \ \ \forall i,j\\
   & b_4\ge 0, b_5\ge 0, q (b_4 + b_5)\le \underline{\sigma}
\end{split}
\end{align}
where $b_4, b_5$ are auxiliary variables and where $q=60$ is the cardinality of the set of delays. Now following the steps in \cite[Appendix B]{xie2022design}, \eqref{measure1}-\eqref{norm2} can be recast as the following
\begin{align}\label{new_constraints}
\begin{split}
      &  [\tilde{T}\tilde{A}_{ii}\tilde{T}^{-1}]_s\le -b_1\I_3, \ \ \forall i\\
       & \left[\begin{matrix}b_3 \I_3 & \frac{\eta_{j}}{\eta_i}(\tilde{T}\tilde{A}_{i{j}}\tilde{T}^{-1})\T \\ \frac{\eta_{j}}{\eta_i}\tilde{T}\tilde{A}_{i{j}}\tilde{T}^{-1} & b_3\I_3\end{matrix}\right]\succeq 0, \ \ \forall i,j\\
       & \left[\begin{matrix}b_4 \I_3 & (\tilde{T}(\tilde{B}_k)_{ii}\tilde{T}^{-1})\T \\ \tilde{T}(\tilde{B}_k)_{ii}\tilde{T}^{-1} & b_4\I_3\end{matrix}\right]\succeq 0, \ \ \forall i\\
       & \left[\begin{matrix}b_5 \I_3 & \frac{\eta_{j}}{\eta_i}(\tilde{T}(\tilde{B}_k)_{i{j}}\tilde{T}^{-1})\T \\ \frac{\eta_{j}}{\eta_i}\tilde{T}(\tilde{B}_k)_{i{j}}\tilde{T}^{-1} & b_5\I_3\end{matrix}\right]\succeq 0,\ \ \forall i,j\\
       &  b_1> 0, b_2 \ge 0, b_3\ge 0, b_4\ge 0, b_5\ge 0, 2b_3\le b_2\\
       &  -b_1+b_2\le -\bar{\sigma}, 60 (b_4 + b_5)\le \underline{\sigma}
        \end{split}
\end{align}
We then solve the optimisation problem \eqref{equ: optimisation} replacing \eqref{measure}-\eqref{norm} with \eqref{new_constraints} for a grid of $\alpha$ and $\beta$ and for the fixed $\eta$\footnote{The code solving the optimisation problem and the data, i.e. $\eta$, can be found in \url{https://tinyurl.com/yc3frafb}.}. The optimal solution for this problem is: $k_0=0.7445, k_1=1.3399, k_2=0.5052,k_1^{(\tau)}= 0.00057, k_1^{(\tau)}=0.00076, k_2^{(\tau)}=0.00048$, with corresponding $\alpha=-0.5, \beta=-1$.}

\subsection{Simulation}
\begin{figure}
\centering
\includegraphics[width=0.95\columnwidth]{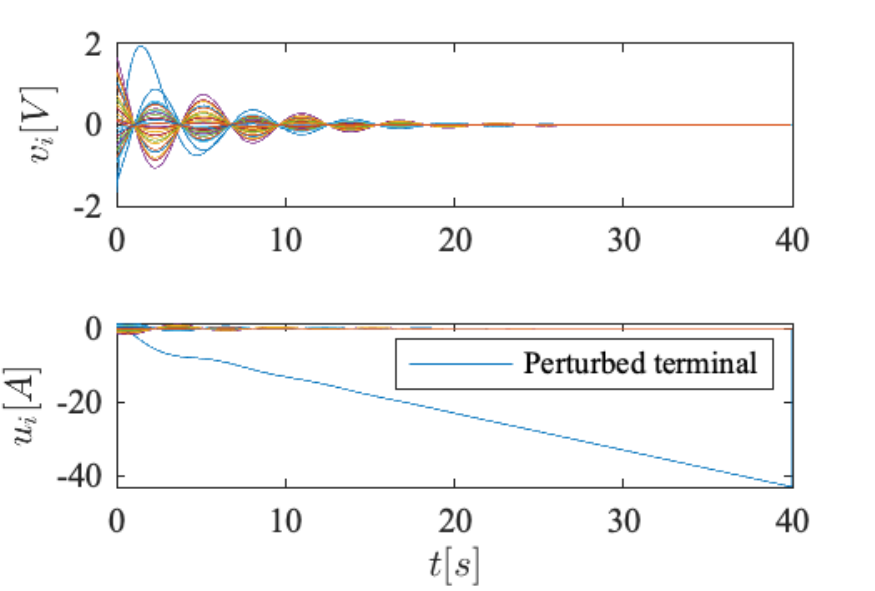}
\caption{Top panel: Voltage deviation from nominal voltage; Bottom panel: Injected control current. }
\label{fig: voltage}
\end{figure}

We consider a disturbance acting on a random terminal, say terminal $1$, which is $d_1(t)=3+t+e^{-0.2t}\sin t$. The heterogeneous delays are selected as $\tau_{ij}(t)=\tau_{i}(t)=0.1+0.1\sin(t+i)$ seconds. The grid is assumed to be initiated at $v_i(0) \sim \sN(0,1), \forall i$. Figure \ref{fig: voltage} (top panel) illustrates the voltage deviation of all the terminals from the nominal value. It shows that all the deviations finally reduce to $0$ including the perturbed terminal $1$. This is in accordance with the theoretical prediction as the designed controller injected a ramp current to compensate for the ramp component in the disturbance, as illustrated in the bottom panel.

\section{Conclusions and future work}
We considered the problem of designing distributed multiplex integral control protocols for nonlinear networks affected by delays and disturbances. The designed control protocol, delivered via multiplex architecture and fulfilling certain conditions leveraging non-Euclidean contraction theory, is able to: (i) reject polynomial disturbances; (ii) achieve Input-to-State Stability for nonlinear networks affected by heterogeneous time-varying delays. We validated the results by considering the problem of controlling an MTDC grid and simulations confirmed the effectiveness of the results. {Apart from devising approaches to systematically determine the {\em right} coordinate transformation matrices $T$ and the weighting matrix $[\eta]^{-1}$, see e.g. \cite{centorrino2023euclidean}, we would like to tackle the challenge of designing protocols that guarantee the stronger scalability property studied in \cite{9353260,xie2022design} for the network systems with heterogeneous delays considered in this paper.}




\begin{thebibliography}{10}

\bibitem{silva2021string}
G.~F. Silva, A.~Donaire, A.~McFadyen, and J.~J. Ford, ``String stable integral
  control design for vehicle platoons with disturbances,'' {\em Automatica},
  vol.~127, p.~109542, 2021.

\bibitem{knorn2014passivity}
S.~Knorn, A.~Donaire, J.~C. Ag{\"u}ero, and R.~H. Middleton, ``Passivity-based
  control for multi-vehicle systems subject to string constraints,'' {\em
  Automatica}, vol.~50, no.~12, pp.~3224--3230, 2014.

\bibitem{lombana2016multiplex}
D.~A.~B. Lombana and M.~Di~Bernardo, ``Multiplex {PI} control for consensus in
  networks of heterogeneous linear agents,'' {\em Automatica}, vol.~67,
  pp.~310--320, 2016.

\bibitem{xie2022design}
S.~Xie and G.~Russo, ``On the design of integral multiplex control protocols
  for nonlinear network systems with delays,'' {\em arXiv preprint
  arXiv:2206.03535 {(submitted to Automatica)}}, 2022.

\bibitem{lohmiller1998contraction}
W.~Lohmiller and J.-J.~E. Slotine, ``On contraction analysis for non-linear
  systems,'' {\em Automatica}, vol.~34, no.~6, pp.~683--696, 1998.

\bibitem{jafarpour2022robust}
S.~Jafarpour, A.~Davydov, M.~Abate, F.~Bullo, and S.~Coogan, ``Robust training
  and verification of implicit neural networks: A non-{Euclidean} contractive
  approach,'' {\em arXiv preprint arXiv:2208.03889}, 2022.

\bibitem{davydov2022non}
A.~Davydov, S.~Jafarpour, and F.~Bullo, ``Non-{Euclidean} contraction theory
  for robust nonlinear stability,'' {\em IEEE Transactions on Automatic
  Control}, 2022.

\bibitem{russo2010global}
G.~Russo, M.~Di~Bernardo, and E.~D. Sontag, ``Global entrainment of
  transcriptional systems to periodic inputs,'' {\em PLoS computational
  biology}, vol.~6, no.~4, p.~e1000739, 2010.

\bibitem{aminzare2014contraction}
Z.~Aminzare and E.~D. Sontag, ``Contraction methods for nonlinear systems: A
  brief introduction and some open problems,'' in {\em 53rd IEEE Conference on
  Decision and Control}, pp.~3835--3847, IEEE, 2014.

\bibitem{di2016convergence}
M.~di~Bernardo, D.~Fiore, G.~Russo, and F.~Scafuti, ``Convergence, consensus
  and synchronization of complex networks via contraction theory,'' {\em
  Complex Systems and Networks}, pp.~313--339, 2016.

\bibitem{tsukamoto2021contraction}
H.~Tsukamoto, S.-J. Chung, and J.-J.~E. Slotine, ``Contraction theory for
  nonlinear stability analysis and learning-based control: A tutorial
  overview,'' {\em Annual Reviews in Control}, vol.~52, pp.~135--169, 2021.

\bibitem{shiromoto2018distributed}
H.~S. Shiromoto, M.~Revay, and I.~R. Manchester, ``Distributed nonlinear
  control design using separable control contraction metrics,'' {\em IEEE
  Transactions on Control of Network Systems}, vol.~6, no.~4, pp.~1281--1290,
  2018.

\bibitem{revay2021recurrent}
M.~Revay, R.~Wang, and I.~R. Manchester, ``Recurrent equilibrium networks:
  Flexible dynamic models with guaranteed stability and robustness,'' {\em
  arXiv preprint arXiv:2104.05942}, 2021.

\bibitem{wang2006contraction}
W.~Wang and J.-J. Slotine, ``Contraction analysis of time-delayed
  communications and group cooperation,'' {\em IEEE Transactions on Automatic
  Control}, vol.~51, no.~4, pp.~712--717, 2006.

\bibitem{5717887}
G.~Russo, M.~di~Bernardo, and E.~D. Sontag, ``Stability of networked systems: A
  multi-scale approach using contraction,'' in {\em 49th IEEE Conference on
  Decision and Control (CDC)}, pp.~6559--6564, 2010.

\bibitem{FB-CTDS}
F.~Bullo, {\em Contraction Theory for Dynamical Systems}.
\newblock Kindle Direct Publishing, {1.0}~ed., 2022.

\bibitem{wen2008generalized}
L.~Wen, Y.~Yu, and W.~Wang, ``Generalized {Halanay} inequalities for
  dissipativity of {Volterra} functional differential equations,'' {\em Journal
  of Mathematical Analysis and Applications}, vol.~347, no.~1, pp.~169--178,
  2008.

\bibitem{siljak2011decentralized}
D.~D. Siljak, {\em Decentralized control of complex systems}.
\newblock Courier Corporation, 2011.

\bibitem{desoer1972measure}
C.~Desoer and H.~Haneda, ``The measure of a matrix as a tool to analyze
  computer algorithms for circuit analysis,'' {\em IEEE Transactions on Circuit
  Theory}, vol.~19, no.~5, pp.~480--486, 1972.

\bibitem{axelson2022online}
M.~Axelson-Fisk and S.~Knorn, ``Online distributed design for control cost
  reduction,'' {\em Automatica}, vol.~141, p.~110312, 2022.

\bibitem{axelson2022graph}
M.~Axelson-Fisk and S.~Knorn, ``A graph theoretic approach to ensure a scalable
  performance measure in multi-agent systems with {MIMO} agents,'' {\em
  IFAC-PapersOnLine}, vol.~55, no.~13, pp.~103--108, 2022.

\bibitem{andreasson2014distributed}
M.~Andreasson, M.~Nazari, D.~V. Dimarogonas, H.~Sandberg, K.~H. Johansson, and
  M.~Ghandhari, ``Distributed voltage and current control of multi-terminal
  high-voltage direct current transmission systems,'' {\em IFAC Proceedings
  Volumes}, vol.~47, no.~3, pp.~11910--11916, 2014.

\bibitem{sneath2014fault}
J.~Sneath and A.~D. Rajapakse, ``Fault detection and interruption in an earthed
  {HVDC} grid using {ROCOV} and hybrid {DC} breakers,'' {\em IEEE Transactions
  on Power Delivery}, vol.~31, no.~3, pp.~973--981, 2014.

\bibitem{centorrino2023euclidean}
V.~Centorrino, A.~Gokhale, A.~Davydov, G.~Russo, and F.~Bullo, ``Euclidean
  contractivity of neural networks with symmetric weights,'' {\em arXiv
  preprint arXiv:2302.13452}, 2023.

\bibitem{9353260}
S.~Xie, G.~Russo, and R.~H. Middleton, ``Scalability in nonlinear network
  systems affected by delays and disturbances,'' {\em IEEE Transactions on
  Control of Network Systems}, vol.~8, no.~3, pp.~1128--1138, 2021.

\end{thebibliography}
\bibliographystyle{ieeetr}
\end{document}